\documentclass[letterpaper,11pt,psamsfonts,reqno]{amsart}

\usepackage{amssymb}
\usepackage{amsxtra}
\usepackage[mathscr]{eucal}

\font\sb = cmbx8 scaled \magstep0

\numberwithin{equation}{section}

\newtheorem{theorem}{Theorem}[section]
\newtheorem{proposition}[theorem]{Proposition}
\newtheorem{lemma}[theorem]{Lemma}
\newtheorem{corollary}[theorem]{Corollary}

\theoremstyle{definition}

\newcommand{\Ad}{\operatorname{Ad}}
\newcommand{\fg}{\mathfrak g}

\newcommand{\h}{\textsf{\textup{h}}}

\newcommand{\Id}{\operatorname{Id}}

\newcommand{\GL}{\operatorname{GL}}
\newcommand{\SL}{\operatorname{SL}}

\newcommand{\htop}{\operatorname{h_{top}}}
\newcommand{\dimbox}{\operatorname{dim_{box}}}

\newcommand{\ba}{\mathbf{a}}
\newcommand{\bb}{\mathbf{b}}

\newcommand{\bm}{\mathbf{m}}
\newcommand{\bn}{\mathbf{n}}

\newcommand{\bs}{\mathbf{s}}
\newcommand{\bt}{\mathbf{t}}

\newcommand{\bv}{\mathbf{v}}

\newcommand{\RR}{\mathbb{R}}

\newcommand{\NN}{\mathbb{N}}
\newcommand{\QQ}{\mathbb Q}
\newcommand{\TT}{\mathbb{T}}
\newcommand{\ZZ}{\mathbb{Z}}

\newcommand{\cD}{\mathcal{D}}

\newcommand{\hd}{Hausdorff dimension}

\newcommand{\ggm}{G/\Gamma}
\newcommand{\nz}{\smallsetminus\{0\}}
\newcommand {\equ}[1]     {\eqref{#1}}
\newcommand\df{\stackrel{\mathrm{def}}{=}}
\newcommand {\ignore}[1]  {}
\newif\ifdraft\drafttrue

\long\def\comd#1{\ifdraft{(\sb #1)}\else\ignorespaces\fi}

\newsymbol\emptyphi 203F

\begin{document}

\title{Measure rigidity and $p$-adic Littlewood-type problems}

\author{Manfred Einsiedler and Dmitry Kleinbock}

\address{Manfred Einsiedler, Department of Mathematics,
Princeton University
Fine Hall, Washington Road,
Princeton NJ 08544}
\email{meinsied@princeton.edu}

\address{Dmitry Kleinbock, Department of Mathematics,
MS 050, Brandeis University, P.O. Box 9110, Waltham, MA 02454}
\email{kleinboc@brandeis.edu}

\date{March  2005}




\begin{abstract} The paper investigates various $p$-adic versions of 
Littlewood's
conjecture, generalizing a set-up considered recently by de~Mathan and
Teuli\'e. In many cases it is shown that the sets of exceptions to 
these conjectures
have \hd\ zero. The proof follows the measure ridigity approach of
Einsiedler, Katok and Lindenstrauss.\end{abstract}

\maketitle

\section{Introduction}\label{sec: intro}

\subsection{The conjecture and our main theorems}
  Every real number $u $ satisfies 
\[
|qu-{q_0}
|<\frac{1}{q}\mbox{ for infinitely many pairs }
(q,q_0)\in\ZZ^2
\]
(where $\,q\neq 0$ is understood). Here the exponent of $q $ is
sharp in the sense that for certain ({\em badly approximable}) numbers,
e.g.\ for irrational elements of
a quadratic number field, there exists a constant $\epsilon>0$
such that
\begin{equation*}\label{e:BA}
|qu-{q_0}
|\geq\frac{\epsilon}{q}
\quad\mbox{ for all pairs }(q,q_0)\in\ZZ^2\,;
\end{equation*}
or, equivalently,
  \[
\liminf_{q\rightarrow\infty }|q|\cdot\langle qu\rangle>0\,,
\]
where $\langle z\rangle\df\min_{q_0\in\ZZ}|z-q_0|$ denotes the
distance to the nearest integer.

\medskip

In the following let $p$ be a prime number. B.~de~Mathan and
O.~Teuli\'e \cite{Mathan-Teulie} conjectured that for every
$u\in\RR$ and $\epsilon>0$
\begin{equation*}\label{e:MT}
|qu-{q_0}
|<\frac{\epsilon}{q|q|_p}=\frac{\epsilon}{ q'}
   \quad \mbox{ for infinitely many pairs }(q,q_0)\in\ZZ^2,
\end{equation*}
where $q=q'p^k$ for some $k\geq 0$, $q'$ is coprime to $p $, and
$|q|_p={1}/{p^k}$ denotes the $p $-adic norm. Equivalently one can
ask whether
\begin{equation}\label{e:MTm}
  \liminf_{q\rightarrow\infty }|q|\cdot|q|_p\cdot\langle qu\rangle=0\,,
\end{equation}
or (see \cite[Lemme 1.3]{Mathan-Teulie}) whether the digits in the
continued fraction expansions for $u,\, pu,\, p^2u,\ldots$ do not all
allow a uniform bound. Moreover, 
de~Mathan and 
Teuli\'e
showed that this 
property holds for elements of a quadratic number field.

We show in the present paper that the above property holds except
possibly for a set of \hd\ zero:

\begin{theorem}\label{theorem: dmt}
The set of $u\in\RR$
which do not satisfy \equ{e:MTm}
  has Hausdorff dimension  zero; in fact, it is a
  countable union of sets with box dimension zero.
\end{theorem}

The above theorem is derived from the following result:

\begin{theorem}\label{theorem: Main}
For any  $v\in\QQ_p$,
the set of $u\in\RR$
which do not satisfy
\begin{equation}\label{g-conjecture}
  \liminf_{
q\rightarrow\infty,\,q_0\rightarrow\infty
}  |q|\cdot|qu-q_0|\cdot|qv-q_0|_p=0
\end{equation}
  has Hausdorff dimension  zero; in fact, it is a
  countable union of sets with box dimension zero.
\end{theorem}



It is easy to see that
Theorem \ref{theorem: Main}
implies Theorem \ref{theorem: dmt}. Indeed, if $u\ne 0$ satisfies
\equ{g-conjecture}
with $v = 0$, then, dividing by $u$ and interchanging $q$ and $q_0$, one gets
\begin{equation}\label{g-workout}
\liminf_{q\rightarrow\infty,\,q_0\rightarrow\infty} 
|q_0|\cdot|\tfrac{q}u - q_0|\cdot|q|_p=0
\,.
\end{equation}
We claim that when $q_0 \ne 0$ and  the above product is small enough,
$\left|\frac{q}{u}\right|$ must be smaller than $2|q_0|$: otherwise
$$|q_0|\cdot|\tfrac{q}u - q_0|\cdot|q|_p
\ge |\tfrac{q}u - q_0|\cdot|q|_p
\ge\tfrac{1}{2}\left|\tfrac{q}{u}\right|\cdot|q|_p\ge \tfrac{1}{2|u|}\,,$$
contradicting \equ{g-workout}. Hence  \equ{e:MTm} holds with $u$ replaced
by $1/u$.
\medskip

An assertion stronger than Theorem~\ref{theorem: Main}, namely that
the set of pairs $\,(u,v)\in\RR\times\QQ_p$
which do not satisfy \equ{g-conjecture} is a
  countable union of sets with box dimension zero, appears to be within
reach, and will be the subject of a forthcoming paper, see \S 
\ref{forthcoming}.
It seems natural to conjecture, generalizing the conjecture
of de~Mathan and
Teuli\'e, that \equ{g-conjecture} holds
for any $(u,v)\in\RR\times\QQ_p$ except for $u = 0$, $v = 0$.

\subsection{Connection to Littlewood's conjecture}\label{sec: conn to L}

The  above conjecture, as well as its weaker version
due to  de~Mathan and
Teuli\'e,  
look very similar to Littlewood's conjecture (c.\ 1930) which states
\begin{equation}\label{e:Lwd}
  \liminf_{q\rightarrow\infty}q\langle qu\rangle\langle qv\rangle=0
\end{equation}
for every $(u,v)\in\RR^2$.
However, this resemblance goes deeper.
For example, the result by 
de~Mathan and 
Teuli\'e on quadratic
numbers \cite{Mathan-Teulie} is an analogue of a result by J.\ W.\
S.\ Cassels and H.\ P.\ F.\ Swinnerton-Dyer
\cite{Cassels-Swinnerton-Dyer} that showed \eqref{e:Lwd} for
elements of the same cubic number field (while the quadratic case
is easy in the case of \eqref{e:Lwd}).

More recently A.~Katok, E.~Lindenstrauss and the first named
author \cite{Einsiedler-Katok-Lindenstrauss} used dynamics on
$\SL(3,\RR)/\SL(3,\ZZ)$ to show that Littlewood's conjecture fails
at most on a set of Hausdorff dimension zero. Here the action is the
$\RR^2$-flow realized as left multiplication by positive diagonal
matrices, i.e.\ $\bt=(t_1,t_2)\in\RR^2$ acts as
$$
  \alpha^\bt (x)=
  \begin{pmatrix}e^{-(t_1+t_2)}&&\\&e^{t_1}&\\&&e^{t_2}
   \end{pmatrix}x \quad \mbox{ for every }x\in\SL(3,\RR)/\SL(3,\ZZ).
$$
The proof in \cite{Einsiedler-Katok-Lindenstrauss} can be roughly
divided into three main steps.
\begin{enumerate}
  \item[(i)]  Linking the Diophantine
  conjecture with the dynamics of $\alpha$.
  \item[(ii)] Studying (and, if possible, classifying) 
$\alpha$-invariant and ergodic
  probability measures on  $\SL(3,\RR)/\SL(3,\ZZ)$.
  \item[(iii)] Linking the box dimension of an exceptional set
with an $\alpha$-invariant measure of positive entropy.

\end{enumerate}

Step (i) uses at its heart the following characterization, see
\cite[Sect.~2]{Margulis-Oppenheim-conjecture},
  \cite[Sect.~30.3]{Starkov-00}, or
  \cite[Prop.~12.1]{Einsiedler-Katok-Lindenstrauss}: $(u,v)\in\RR^2$ fails
Littlewood's conjecture \eqref{e:Lwd} if and only if
\[
  \begin{pmatrix}1&&\\u&1&\\v&&1\end{pmatrix}\SL(3,\ZZ)\mbox{ has bounded
   orbit
    under the semigroup }\{\alpha^\bt:\bt\in\RR_+^2\}.
\]
In principle, this could be used to solve Littlewood's conjecture
if the dynamics of $\alpha $ were completely understood. Margulis
\cite{Margulis-conjectures}
conjectured that there are very few $\alpha$-invariant measures as
in (ii) -- a phenomenon called measure rigidity. However, the only
currently available method to study these measures works
exclusively for measures satisfying an additional assumption,
namely for measures $\mu$
such
that for some element
$\alpha^\bt$ of the action its measure-theoretic entropy
$\h_\mu(\alpha^\bt)$ is positive. The complete answer for positive
entropy measures was given in
\cite[Thm.~1.3]{Einsiedler-Katok-Lindenstrauss} where the two
separately developed methods from \cite{Einsiedler-Katok-02} and
\cite{Lindenstrauss-Quantum} were combined. In  step
(iii) the positive entropy assumption is then translated to the
statement about the Hausdorff dimension.

\subsection{Method of proof}

For the proof of Theorem \ref{theorem: Main} 
we follow in
principle the same scheme as described in \S \ref{sec: conn to L}, 
but use dynamics on
\begin{equation}\label{e:X}
  X=G/\Gamma,\mbox{ where }G = \SL(2,\RR)\times\SL(2,\QQ_p),\  \Gamma = \imath
\left(\SL\big(2,\ZZ[\tfrac{1}{p}]\big)\right)
\end{equation}
instead.  Here 
$\imath:\SL\big(2,\ZZ[\tfrac{1}{p}]\big)\rightarrow
G$ is the diagonal embedding, i.e.\
$\imath(A)=(A,A)$ for any $A\in \SL\big(2,\ZZ[\tfrac{1}{p}]\big)
$. The action considered is the $(\RR\times\ZZ)$-action defined by
left multiplication of $x\in X$ by
\begin{equation}\label{e:alpha}
  \alpha^{(t,n)} \df \left(\begin{pmatrix}e^{-t}&\\&e^t\end{pmatrix},
\begin{pmatrix}p^{n}&\\&p^{-n}\end{pmatrix}\right)\mbox{ for 
}(t,n)\in\RR\times\ZZ.
\end{equation}
We prove in  \S\ref{dynamical} (see Proposition \ref{prop: dynamical 
g}) that  $u\in\RR\nz$ and $v\in\ZZ_p$
  fail
\eqref{g-conjecture} if and only if
\begin{equation}\label{e:x} 
x_{u,v}\df\left(\begin{pmatrix}1&\\u&1\end{pmatrix},\begin{pmatrix}1&\\v&1\end{pmatrix}\right)\Gamma
\end{equation}
  has bounded
   orbit
    under the semigroup $\alpha^{C} \df \{\alpha^{(t,n)}: (t,n)\in C\}$, where
\begin{equation}\label{e:C}
C  = \{(t,n)\in\RR\times\ZZ:\ n\ge 0,\ t - n\log p\geq 0\}\,.
\end{equation}
This serves as a replacement for  part (i) of the argument of
\cite{Einsiedler-Katok-Lindenstrauss}, and sets the stage for part (ii),
that is,
using measure rigidity methods. Namely, applying a recent result by
E.~Lindenstrauss \cite[Thm.~1.1]{Lindenstrauss-Quantum} (see
Theorem \ref{theorem: Elon}), 
we prove that
any $\alpha$-invariant and ergodic probability
   measure  on $X$ for which the measure-theoretic entropy
(see \S\ref{rigidity}) with respect  to $\alpha^{(1,0)}$ is positive
has to be the
   unique $G$-invariant Haar measure
   on $X$.
Consequently,  no $\alpha$-invariant and ergodic
  probability measure supported on a compact subset $Y$
of $X$ has positive entropy for $\alpha^{(1,0)}$.

 From this, using a generalization
of step (iii) of the proof in \cite{Einsiedler-Katok-Lindenstrauss}
(see Proposition \ref{pro: box dimension}),
we
derive that for any  compact subset $B$ of the  unstable horospherical subgroup
relative to $\alpha^{(1,0)}$ (defined in \S \ref{entropy}), any compact $Y\subset X$,
and any $x\in X$,  the set
   \[
\big\{h\in B:\alpha^C hx\subset Y\big\}
   \]
   has box dimension zero.  And to finish the proof of Theorem 
\ref{theorem: Main}
it remains to observe that $\{x_{u,v}: u\in\RR\}$ is precisely
the orbit of $x_{0,v}$ under the unstable horospherical subgroup
relative to $\alpha^{(1,0)}$.

\medskip

The paper is organized as follows. The aforementioned three parts of 
the argument
are discussed in detail in \S\S \ref{dynamical}--\ref{entropy}
respectively, and
then used in  \S \ref{proof} where we finish the proof of Theorem 
\ref{theorem: Main}.
After that we mention  several  generalizations and open questions. 
In particular
we point out that, while Littlewood's conjecture does not seem to get 
easier to prove in higher dimensions,
the following modification of  the de~Mathan-Teuli\'e conjecture 
follows quite easily from a theorem of Furstenberg 
\cite{Furstenberg67}:

\begin{theorem}\label{t:Furstenberg}
  Let $p_1,p_2$ be two distinct prime numbers. Then for every $u\in\RR$
  \begin{equation}\label{e:MTmF}
  \liminf_{q\rightarrow\infty 
}|q|\cdot|q|_{p_1}\cdot|q|_{p_2}\cdot\langle qu\rangle=0.
\end{equation}
\end{theorem}

This theorem is proved in \S \ref{s:Furstenberg}. Another 
modification of  \eqref{e:MTm},
where the $p$-adic norm is replaced by a ``pseudo-norm''
$$
|q|_{\cD} = \inf\{1/r_n : q\in r_n\ZZ\}\,,
$$
is discussed in \S \ref{s:D}. Here $\cD = ( r_n)_{n\in\ZZ_+}$
is a sequence of positive integers with $r_0 = 1$ and  $r_{n+1}
\in r_n\ZZ$ for each $n$; the choice $r_n = p^n$ gives the usual $p$-adic norm.
De~Mathan and Teuli\'e proved in \cite{Mathan-Teulie} that every 
quadratic irrational
$u$ satisfies
$ \liminf_{q\rightarrow\infty }|q|\cdot|q|_{\cD}\cdot\langle qu\rangle=0$
for any $\cD$ as above with uniformly bounded ratios $r_{n+1}
/r_n$. In this paper we explain how a modification of our method
yields a $\cD$-adic analogue of
  Theorem \ref{theorem: Main} in the case
$\cD = ( a^n)_{n\in\NN}$, where $a$ is not necessarily prime.

\medskip

{\bf Acknowledgements:}
The authors are grateful to University of Washington and Brandeis University
for arranging their mutual visits, and to
the Max Planck Institute for its hospitality during July 2004. Thanks 
are also due to Yann Bugeaud, 
Elon Lindenstrauss and George Tomanov for helpful discussions. This 
research was partially supported by NSF grants 0400587, 0509350, 
0239463,  and 0400587.

\ignore{In fact  Furstenberg \cite{Furstenberg67} showed that a 
closed $\times p_1,\times p_2$-invariant sets is either a finite set 
of rational points or is everything. The above theorem follows since 
for any $\delta>0$ the set
\[
  F_\delta=\bigl\{u\in\RR: 
\sup_{q\in\NN}|q|\cdot|q|_{p_1}\cdot|q|_{p_2}\cdot\langle 
qu\rangle\geq\delta\bigr\}/\ZZ\subset\RR/\ZZ
\]
is closed, $\times p_1,\times p_2$-invariant (see \S \ref{s:Furstenberg}),
and does not contain any rational points.}

\section{Dynamical characterization of
conditions
\equ{e:MTm} and \equ{g-conjecture}}
\label{dynamical}

\subsection{Lattices in $\RR^2\times\QQ_p^2$ and Mahler's
criterion}\label{sec: Mahler}
   In the following we assume that the metric (and the topology) of
$X = G/\Gamma$ as in (\ref{e:X}) is induced by a
right-invariant metric on 
$G$.

It is well known that
$\ZZ[\tfrac{1}{p}]^2$ diagonally imbedded in $\RR^2\times\QQ_p^2$ is a lattice,
that is, it is  discrete
and the quotient $(\RR^2\times\QQ_p^2)/\ZZ[\tfrac{1}{p}]^2$ has finite volume
(in fact, equal to $1$). For any $g  \in G$,
  $g\ZZ[\tfrac{1}{p}]^2$ is also a lattice in $\RR^2\times\QQ_p^2$
of covolume $1$. In
other words, $g\ZZ[\tfrac{1}{p}]^2$ consists of all elements of the form
\[
\left(A\begin{pmatrix}q\\q_0\end{pmatrix},
B\begin{pmatrix}q\\q_0\end{pmatrix}\right)\mbox{ for
}q,q_0\in\ZZ[\tfrac{1}{p}]\,,
\]
where $A\in\SL(2,\RR)$  and
$B\in\SL(2,\QQ_p)$.
That is, it is
a free $\ZZ[\tfrac{1}{p}]$-module generated by
$(\ba_1,\bb_1), (\ba_2,\bb_2)$, where $A$ has column vectors 
$\ba_1,\ba_2\in\RR^2$ and
$B$ has column vectors $\bb_1,\bb_2\in\QQ_p^2$.

It is easy to see that $\Gamma$ is exactly the stablilizer of 
$\ZZ[\tfrac{1}{p}]^2$
under the action described above, and therefore $X = \ggm$ can be identified
with the $G$-orbit of $\ZZ[\tfrac{1}{p}]^2$. In other words,
points $x=(A,B)\Gamma\in X$ are identified with
unimodular
lattices $\Lambda_x$ in $\RR^2\times\QQ_p^2$ that are
generated, as above, by the column vectors of $A$ and $B$.

\ignore{So $\Lambda$ is
uniquely determined by specifying $(A,B)$, but $(A,B)$ is not
uniquely determined by $\Lambda$. In fact, $(A,B)\imath(C)$ for
$C\in\GL(2,\ZZ[\frac{1}{p}])$ are all the pairs of matrices that
describe the same lattice $\Lambda$.
Therefore, the space of
lattices in $\RR^2\times\QQ_p^2$ is identified with
$$\big(\GL(2,\RR)\times\GL(2,\QQ_p)\big)/\imath\big(\GL(2,\ZZ[{1}/{p}])\big)\,.$$

Thus we can identify points $x=(A,B)\Gamma\in X$ with
(unimodular)
\comd{don't we need to explain unimodular somehow?}
lattices $\Lambda_x$ in $\RR^2\times\QQ_p^2$ that are
generated, as above, by the column vectors of $A$ and $B$.
Note that this way $\Gamma\in X$ is identified with  ${\ZZ}[1/p]^{2}$
diagonally imbedded in $\RR^2\times\QQ_p^2$.}

\medskip

This correlation between points in $X$ and certain lattices allows
a convenient description of compact subsets of $X$ (see
\cite[Thm.\ 7.10]{Kleinbock-Tomanov} for a more general formulation of
the $S$-adic Mahler's criterion).

\begin{theorem}[Mahler's criterion for $X$]\label{thm: MC}
  A subset $L\subset X$ has compact
closure if and only if there exists some $\delta>0$ so that
\[
  L\subset K_\delta\df\Bigl\{x\in X: \Lambda_x\cap
  B_\delta^{\RR^2\times\QQ_p^2}(0)=\{0\}\Bigr\}.
\]
\end{theorem}


\begin{proof} The implication ($\Longrightarrow$) is trivial (and not needed
for our purposes).
For the converse, we need to show that $K_\delta$ is
compact for every positive $\delta$. Let us denote 
$\SL(2,\RR)\times\SL(2,\ZZ_p)$ by $G'$,
and observe that $G = G'\Gamma$ due to  the strong
approximation theorem 
        \cite{Kn}; in other words,
   every $g \in
G$ can be represented as $g =
g_{f}g_l$, where $g_{f} \in  G'$ and
$g_{l} \in
\Gamma$. Hence \begin{equation*}
X\cong G'/(G'\cap \Gamma) = G'/\SL(2,\ZZ)\cong\{g\ZZ^2 : g\in G'\}\,.
\end{equation*}
Furthermore, one has
\begin{equation*}
\begin{aligned}
g\ZZ[\tfrac{1}{p}]^2 \cap  (\RR\times\ZZ_p)^{2} &=
g_{f}\ZZ[\tfrac{1}{p}]^2 \cap  (\RR^{2}\times\ZZ_p^{2})\\& =
g_{f}\left(\ZZ[\tfrac{1}{p}]^2 \cap  (\RR^{2}\times\ZZ_p^{2})\right)
= g_{f}{\ZZ}^2\,.%
\end{aligned}
\end{equation*}
Thus one can identify $K_\delta$ with
\begin{equation}\label{e:mahler}
\left\{g\ZZ^2 : g\in G',\ g {\ZZ}^2 \cap
  B_\delta^{\RR^2\times\ZZ_p^2}(0)=\{0\}\right\}\,.
\end{equation}
To finish the proof it remains to observe that the projection
$G'\to \SL(2,\RR)$ induces a surjection $G/\Gamma \to \SL(2,\RR)/\SL(2,\ZZ)$
with compact fibers, and the set \equ{e:mahler} is contained in the preimage
of
\begin{equation*}
\left\{g\ZZ^2 : g\in \SL(2,\RR),\ g {\ZZ}^2 \cap
  B_\delta^{\RR^2}(0)=\{0\}\right\}\,,
\end{equation*}
which is compact due to to the original Mahler's Criterion
(see \cite[Corollary 10.9]{R}).
\end{proof}

\subsection{Theorem \ref{theorem: Main} and cone orbits}

A {\em cone $C$} of $\RR\times \ZZ$ is a subset defined by two
inequalities $(t,n)\cdot \bv_1\geq 0$ and $(t,n)\cdot \bv_2\geq 0$
where $\bv_1,\bv_2\in\RR^2$ are linearly independent vectors. The
{\em cone orbit} of a point $x\in X$ is given by
\[
  \alpha^Cx=\{\alpha^{(t,n)}x:(t,n)\in C\}\,,
\]
where $\alpha$ is as defined in \equ{e:alpha}.

\begin{proposition}\label{prop: dynamical g}
  Let $C$ be as in \equ{e:C}, and let $u\in \RR$ be nonzero. Then
  $(u,v)$ satisfies
\eqref{g-conjecture}
if and only if
the $C$-orbit of
$
   x_{u,v}$
as in  \equ{e:x}
  is unbounded.
\end{proposition}

We will only prove the `if' part, since this is
the only direction that will be needed for Theorems \ref{theorem:
Main}. The converse is not difficult either.

\begin{proof}
  Suppose $x_{u,v}$ has unbounded $C$-orbit. By Mahler's criterion
(Theorem \ref{thm: MC})
  there exists for every $\delta>0$ a pair $(t,n)$ with $n\ge 0$
  and $e^tp^{-n}\geq 1$ such that $\Lambda_{\alpha^{(t,n)}x_{u,v}}$
  contains a nonzero element in $B_\delta^{\RR^2\times\QQ_p^2}(0)$.
  Clearly, $\Lambda_{\alpha^{(t,n)}x_{u,v}}$ is generated by
  \[
   \left(\begin{pmatrix}e^{-t}\\ e^tu\end{pmatrix},
   \begin{pmatrix}p^{n}\\ p^{-n}v \end{pmatrix}\right)\mbox{ and }
   \left(\begin{pmatrix}0\\ e^t\end{pmatrix},
   \begin{pmatrix}0\\p^{-n}\end{pmatrix}\right).
  \]
  However, since $\Lambda_{\alpha^{(t,n)}x_{u,v}}$ is a
  $\ZZ[\frac{1}{p}]$-module, the vectors
  \[
   \left(\begin{pmatrix}e^{-t}p^{-n}\\ e^tp^{-n}u\end{pmatrix},
   \begin{pmatrix}1\\ p^{-2n}v \end{pmatrix}\right)\mbox{ and }
   \left(\begin{pmatrix}0\\ e^tp^{-n}\end{pmatrix},
   \begin{pmatrix}0\\p^{-2n}\end{pmatrix}\right)
  \]
  are also generators. Therefore, there exists some nonzero
  $(q,q_0)\in\ZZ[\frac{1}{p}]^2$ such that
  \[
   \left(\begin{pmatrix}e^{-t}p^{-n}q\\ e^tp^{-n}(qu-q_0)\end{pmatrix},
   \begin{pmatrix}q\\p^{-2n}(qv-q_0)\end{pmatrix}\right)\in\RR^2\times\QQ_p^2
  \]
  is $\delta$-small. In particular,
$
|q|_p$ is less than $\delta$,
which (assuming $\delta \le 1$) implies that $q\in\ZZ$. Since
  $n\geq 0$ and $v\in\ZZ_p$, the inequalty
\begin{equation}
\label{e-q0}
p^{2n}|(qv-q_0)|_p = |p^{-2n}(qv-q_0)|_p <\delta
\end{equation}
shows that  $q_0\in\ZZ$ as well.

Now consider two cases. If
$u\in\QQ$, that is, $qu - q_0 = 0$ for some $q,q_0\in\ZZ\nz$,
\eqref{g-conjecture} is
obviously satisfied (simply replace
$q,q_0$ in \eqref{g-conjecture}
by $kq,kq_0$, $k\in\NN$).
  If
$u\notin\QQ$, we use the assumption $e^tp^{-n}\geq 1$ to conclude that
the inequality
\begin{equation}
\label{e-q}
|e^tp^{-n}(qu-q_0)| <\delta
\end{equation}
can only hold for $|q|$ and $|q_0|$ not less than some $Q = Q(u,\delta)$, where
the latter quantity tends to infinity as $\delta \to 0$. By taking the
product of the
  inequality $|e^{-t}p^{-n}q| <\delta$ with \equ{e-q} and \equ{e-q0},
we arrive at
  $|q|\cdot|qu-q_0|\cdot|qv-q_0|_p<\delta^3$. Since this holds for any $\delta$,
  we conclude that \eqref{g-conjecture} holds.
\end{proof}

\subsection{The de Mathan-Teuli\'e Conjecture  and cone orbits} Even though
we were able to reduce Theorem \ref{theorem: dmt} to Theorem 
\ref{theorem: Main}, it is
instructive to observe that condition \equ{e:MTm} can also be characterized
in a dynamical language.

\begin{proposition}\label{prop: dynamical MT}
  Let $C'=\{(t,n)\in\RR\times\ZZ: n\leq 0, t+n\log p\geq 0\}$. Then
  $u\in\RR\nz$ satisfies \eqref{e:MTm} if (and only if)
  $
   x_{u,0}$
  has unbounded $C'$-orbit, i.e.\ $\alpha^{C'}x_{u,0}\subset X$ does not
  have compact closure.
\end{proposition}

\begin{proof}
  We again only prove the if-part. As in the proof of Proposition
  \ref{prop: dynamical g} we find for an arbitrary $\delta>0$ some $(t,n)\in C'$
  such that $\Lambda_{\alpha^{(t,n)}x_{u,0}}$ contains some nonzero
  element of $B_\delta^{\RR^2\times\QQ_p^2}(0)$. We multiply the 
obvious generators
  of $\Lambda_{\alpha^{(t,n)}x_{u,0}}$ by $p^n$ and get the two
  generators
  \[
   \left(\begin{pmatrix}e^{-t}p^{n}\\ e^tp^{n}u\end{pmatrix},
   \begin{pmatrix}p^{2n}\\ 0 \end{pmatrix}\right)\mbox{ and }
   \left(\begin{pmatrix}0\\ e^tp^{n}\end{pmatrix},
   \begin{pmatrix}0\\1\end{pmatrix}\right).
  \]
  Therefore, the small element has the form
  \[
   \left(\begin{pmatrix}e^{-t}p^{n}q\\ e^tp^{n}(qu-q_0)\end{pmatrix},
   \begin{pmatrix}p^{2n}q\\-q_0\end{pmatrix}\right)\in\RR^2\times\QQ_p^2
  \]
  for some $(q,q_0)\in\ZZ[\frac{1}{p}]^2$. Since $n\leq 0$ we
  conclude from the $p$-adic part of that vector that $q,q_0\in\ZZ$,
and in fact
\begin{equation}
\label{e-q0-MT}
p^{-2n}|q|_p  <\delta
\end{equation}
  By definition of $C'$ we have $e^tp^n\geq 1$,
thus, in the case $u\notin\QQ$ (the rational case is again easy)
the inequality
\begin{equation}
\label{e-q-MT}
|e^tp^{n}(qu-q_0)| <\delta
\end{equation}
forces $|q|$ and $|q_0|$ to big enough. By taking the
product of the
  inequality $|e^{-t}p^{n}q| <\delta$ with \equ{e-q0-MT} and \equ{e-q-MT},
we arrive at
  $|q|\cdot|q|_p\cdot|qu-q_0|<\delta^3$.
\end{proof}

\section{Measure rigidity on $X$}
\label{rigidity}

Measure rigidity refers to the phenomenon that certain actions have very few
invariant measures. The case of unipotent actions is well understood due to
Ratner's work \cite{Ratner-measure-rigidity}, and extensions by Ratner
\cite{Ratner-padic} and
Margulis-Tomanov \cite{Margulis-Tomanov}. The case of higher rank 
partially hyperbolic
actions,
as it concerns us here, has seen recently an interesting development
\cite{Kalinin-Katok-99, Kalinin-Spatzier-02, Einsiedler-Katok-02, Lindenstrauss-Quantum, 
Einsiedler-Katok-Lindenstrauss}.
  However  here our understanding is not complete, and positive 
entropy is, so far, a crucial
  additional assumption. Even so it is possible to apply these 
results, in particular to
  number theory.
\ignore{For Theorem \ref{theorem: Main}
  we will use a result by E.~Lindenstrauss \cite{Lindenstrauss-Quantum}, and
for Theorem \ref{theorem: More}
  we will need its generalization, to be stated later.}

Let $H$ be a locally compact metric group
acting on a locally compact metric space $X$. A
measure $\mu$ on $X$ is {\em $H$-invariant} if
$\mu(h^{-1}B)=\mu(B)$ for any measurable $B\subset X$. A measure $\mu$
is {\em $H$-ergodic} if any invariant measurable set, that is any
measurable $B\subset X$ with
$\mu(B\,\Delta\, hB)=0$
for any $h\in H$,
satisfies $\mu(B)\in\{0,1\}$. We will only consider probability
measures. Then the $H$-invariant and ergodic measures are the
extremal points of the convex set of all $H$-invariant measures.
Moreover, every $H$-invariant measure $\nu$ can be expressed as a
generalized convex combination
\begin{equation}\label{ergodic components}
  \nu=\int_R \mu_r\operatorname{d}\!\rho(r)
\end{equation}
of $H$-invariant and ergodic measures $\mu_r$ on $X$ ($r\in R$),
where $(R,\rho)$ is some probability space. This is the {\em
ergodic decomposition} of $\mu$ with respect to $H$, and the
measures $\mu_r$ are the {\em ergodic components} \cite{Varadarajan-63}.

We will need the notion of measure theoretic entropy. Let
$T:X\rightarrow X$ be a measure preserving map on the probability
space $(X,\mu)$. Instead of giving the formal definition
\cite[Ch.~4]{Walters-82} let us just mention that the measure
theoretic entropy $\h_\mu(T)$ is a non-negative (in general
possibly infinite) number that measures the complexity of the
dynamical system defined by $T$ on $(X,\mu)$.

As mentioned in the introduction, the following theorem can be
specialized to fit our needs. To avoid unnecessary complications we 
state a slightly
simplified version.

\begin{theorem}\cite[Thm.\ 1.1]{Lindenstrauss-Quantum}\label{theorem: Elon}
    Let $H = \SL(2,\RR)$, let $L$ be an $S$-arithmetic group, and let 
$G=H \times L$.
Take $\Gamma$ to be a discrete subgroup of $G$ (not necessarily a
lattice) such that $\Gamma\cap L$ is finite. Suppose $\mu$ is a
probability measure on $X = G/\Gamma$ that is
invariant under left multiplication by elements of the diagonal
group $A$ in $\SL(2,\RR)$. Assume furthermore the following two
conditions.
\begin{enumerate}\item All ergodic components of $\mu$ with
respect to the action of $A$ have positive entropy.
\item $\mu$ is
$L$-recurrent, i.e.\ for a measurable $B\subset X$ a.e.\ $x\in B$ 
satisfies for every compact subset
  $C\subset L$ that there exists $\ell\in L\smallsetminus C$ with $\ell x\in B$.
\end{enumerate} Then $\mu$ is a convex
combination of algebraic measures invariant under $H$.
\end{theorem}

For the case $G=\SL(2,\RR)\times\SL(2,\QQ_p)$ and the lattice
$\Gamma=\imath\left(\SL\big(2,\ZZ[\frac{1}{p}]\big)\right)$ we obtain the
following corollary.

\begin{corollary}\label{corollary}
   Let $X=G/\Gamma$ be as in \eqref{e:X}, and let
   $\alpha$ be the $(\RR\times\ZZ)$-action defined in  \eqref{e:alpha}.
   Let $\mu$ be an $\alpha$-invariant and ergodic probability
   measure on $X$. If $\h_\mu(\alpha^{(1,0)})>0$, then $\mu$ is the
   unique $\SL(2,\RR)\times\SL(2,\QQ_p)$-invariant Haar measure
   $\mu=m_X$ of $X$.
\end{corollary}

The following general property of entropy will be useful (see
\cite[Thm.~8.2]{Walters-82} in the case of the ergodic decomposition and
\cite{Jacobs-63} for the general case). If
$T:X\rightarrow X$ is continuous, $\nu$ is $T$-invariant, and
$\nu=\int_R\mu_r\operatorname{d}\!\rho(r)$ is a convex combination
of other $T$-invariant probability measures, then
\begin{equation}\label{e: entropy is affine}
   \h_\nu(T)=\int\h_{\mu_r}(T)\operatorname{d}\!\rho(r).
\end{equation}
In particular, the entropy $\h_\nu(T)$ for a non-ergodic measure
$\nu$ as in \eqref{ergodic components} can be calculated from the
entropy $\h_{\mu_r}(T)$ of the ergodic components $\mu_r$ (and the
measure $\rho$).

\begin{proof}
   We will apply Theorem \ref{theorem: Elon} to $L=\SL(2,\QQ_p)$, and 
$\Gamma$ as above.
   Clearly, $\Gamma\cap L$ is the trivial group. Now let $\mu$ be an
   $\alpha$-invariant and ergodic probability measure on $X$. Then
   it is invariant under $A$. Note however, that $\mu$ might not be
   ergodic under $A$. So by the ergodic decomposition
   $\mu$ possibly decomposes in a generalized
   convex combination $\mu=\int_R \nu_r\operatorname{d}\!\rho$
   of $A$-ergodic measures $\nu_r$, where $(R,\rho)$ is some probability
   space. However, ergodicity of $\mu$ under $\alpha$ implies that
   $\h_{\nu_r}(\alpha^{(1,0)})=\h_\mu(\alpha^{(1,0)})$ for
   $\rho$-a.e.\ $r\in R$. By assumption this entropy is positive, so
   condition (1) in Theorem \ref{theorem: Elon} is satisfied.

   Since $\mu$ is a probability measure and invariant under
   $\alpha^{(0,1)}$, it is also recurrent under $\alpha^{(0,1)}$ by
   Poincar\'e recurrence. Since $\alpha^{(0,1)}$ is left
   multiplication by an element of $L$ (whose powers are not
   contained in a compact subset of $L$), it follows that $\mu$ is
   recurrent under $L$ as required in condition (2).

   Therefore, $\mu$ is a (generalized) convex combination of
   homogeneous measures that are invariant under $\SL(2,\RR)$ by
   Theorem \ref{theorem: Elon}.  The corollary follows from the
    following proposition.
  \end{proof}

  \begin{proposition}
    Let $X=G/\Gamma$ be as in \eqref{e:X}.
    The Haar measure $m_X$
   is the only homogeneous measure $\nu$ on $X$
   that is invariant under $\SL(2,\RR)$.
\end{proposition}

\begin{proof}
  The proposition follows from \cite[Thm.~1]{Tomanov-00} since $\Gamma$ is
  an $S$-arithmetic lattice, and since there doesn't exist a proper finite index
  subgroup of $G$. We also give a short self-contained proof.

   Let $H'\subset G$ be a closed subgroup that contains $\SL(2,\RR)$ such
  that $\nu$ is the Haar measure on a single closed $H'$-orbit $H'x$.
  Let $x=g\Gamma$ for some $g=(g_\infty,g_p)\in G$, then $H'x\simeq H'/\Lambda'$
  where $\Lambda'=g\Gamma g^{-1}\cap H'$ is a lattice in $H'$. 
Clearly, $H''=H'\cap
  \bigl(\SL(2,\RR)\times\SL(2,\ZZ_p)\bigr)$ is an open closed subgroup
  of $H'$. Therefore, $\Lambda''=\Lambda'\cap H''$ is a lattice in $H''$,
  and the factor map $\pi:H''\rightarrow \SL(2,\RR)$ with compact kernel
  $H''\cap\SL(2,\ZZ_p)$ descends to a factor map $\pi:H''/\Lambda'\rightarrow
  \SL(2,\RR)/\pi(\Lambda'')$. It follows that $\pi(\Lambda'')\subset\SL(2,\RR)$
  is a lattice. However, since $\Gamma_g=
  \pi\Big(g\Gamma g^{-1}\cap \big(\SL(2,\RR)\times \SL(2,\ZZ_p)\big)\Big)$ is
  also a lattice (that is commensurable to $g_\infty \SL(2,\ZZ) g_\infty^{-1}$),
  it follows that $\pi(\Lambda'')\subset \Gamma_g$ has finite index.
  Let $$u_N=\begin{pmatrix}1&N\\&1\end{pmatrix}$$ be an integer matrix 
(which we can
  consider as an element of $\SL(2,\RR)$ as well as $\SL(2,\QQ_p)$), and let
  $u_N^T$ denote its transpose.
  Then for sufficiently large $N$ we have $g_\infty u_N g_\infty^{-1},
  g_\infty u_N^T g_\infty^{-1}\in\pi(\Lambda'')$, and so
  by the definition of $\Gamma$ we have 
$g_pu_Ng_p^{-1},g_pu_N^Tg_p^{-1}\in H''$.
  These two elements generate a subgroup whose closure is $g_p K_n g_p^{-1}$
  for some $n$, where
  $$ K_n=\bigl\{ (\Id,A): A\in\SL(2,\ZZ_p),A\equiv\Id\mbox{ (mod }p^n) \bigr\}.
  $$
  We conclude that $H'$ contains $K_m$ for some $m$. Using the definition
  of $\Gamma$ it follows that $H'x=X$ and that the unique $H'$-invariant
  measure $\nu$ on $X$ is the Haar measure of $X$.
\end{proof}

\section{The box dimension for points with cone orbit inside a compact set}
\label{entropy}

In the case of a single continuous map $T:X\rightarrow X$ on a 
compact metric space
it is well known that every orbit is equidistributed
(with respect to $\mu$) if and only if $\mu$ is the only invariant
probability measure on $X$. In the case of a statement like
Theorem \ref{theorem: Elon} we cannot make a statement about every orbit,
we can, however, deduce that the set of points with bounded orbits
has small dimension.

We will give the result in a more general setting than is used in 
this paper. For this
let $G=\prod_{\sigma\in S}G_\sigma$ be a product of real and $p$-adic 
Lie groups
$G_\sigma$,
where $S\subset\{\infty,p:p\in\NN$ is prime$\}$ is a finite set of 
places of $\QQ$,
let $\Gamma\subset G$ be
a discrete subgroup, and let $\alpha:\RR^k\times\ZZ^\ell\rightarrow G$ be
a homomorphism such that the image consists of semisimple elements of $G$.
Then we consider $\alpha$ as an $\RR^k\times\ZZ^\ell$-action on $X= G/\Gamma$.
Let $d(\cdot,\cdot)$ be a right invariant metric on $G$
which also induces a metric
on $X$, again denoted by $d(\cdot,\cdot)$. We will assume that the metric $d$
is such that the locally defined maps $\log$ and $\exp$ between $G$ and its
Lie algebra are Lipshitz. Here the Lie algebra is the product of the real
and $p$-adic Lie algebras $\fg_\sigma$ corresponding to the factors 
$G_\sigma$ for $\sigma\in S$, and
we equip it with metric derived from norms on the real and $p$-adic 
Lie algebras $\fg_\sigma$.

For a fixed $(\bt,\bn)\in \RR^k\times\ZZ^\ell$ we let
$$
U^+=\bigl\{g\in G: \alpha^{j(\bt,\bn)}g\alpha^{-j(\bt,\bn)}\rightarrow e
  \mbox{ for }j\rightarrow -\infty\bigr\}
$$
be the horospherical unstable subgroup corresponding to
$\alpha^{(\bt,\bn)}$, the horospherical stable
subgroup is defined similarly. Clearly, $U^+x$ gives the unstable `manifold'
(in the presence of $p$-adic Lie groups it is not a manifold in the 
usual sense)
for $\alpha^{(\bt,\bn)}$ through $x\in X$.

A subset $C\subset\RR^k\times\ZZ^\ell$ defined by
$C=\{(\bt,\bn): (\bt,\bn)\cdot\bv_i\geq 0\}$, where
$\bv_i\in\RR^{k+\ell}$ are given vectors for $i=1,\ldots,I$,
is  a {\em cone} if the group generated by $C$ is $\RR^k\times\ZZ^\ell$.
We will be interested in the parts of the orbits corresponding to $C$.
Here we allow $I=0$ and $C=\RR^k\times\ZZ^\ell$.

Before we state the result proved in this section, let us recall
the definition of box dimension of a subset $E$ of a metric space
$(Z,d)$. For $\delta>0$ let $s_\delta(E)$ be the maximal number of
points $u_1,u_2,\ldots$ such that $d(u_i,u_j)>\delta$ for any
$i\neq j$, then the {\em (upper) box dimension} is
\[
   \dimbox(E)=\limsup_{\delta\rightarrow 0}\frac{\log
   s_\delta(E)}{|\log\delta|}.
\]

\begin{proposition}\label{pro: box dimension}
   Suppose $\alpha$ and $X=G/\Gamma$ are as above,
   and $C\subset\RR^k\times\ZZ^\ell$ is a cone that contains $(\bt,\bn)$.
  Let $Y\subset X$ be a compact set such that no $\alpha$-invariant and ergodic
  probability measure supported on $Y$ has positive entropy for 
$\alpha^{(\bt,\bn)}$.
  Then for any compact $B\subset U^+$ of the unstable horospherical subgroup
  and any $x\in X$  the set
   \[
    E_{B,x}(Y)=\left\{u\in B:\alpha^C ux\subset Y\right\}
   \]
   has box dimension zero. Moreover, if also $-(\bt,\bn)\in C$, then
the same is true for compact subsets $B\subset U^+U^-$
and the box dimension of
\[
    Y_C=\bigl\{y\in Y:\alpha^Cy\subset Y\bigr\}
\]
is at most $\dimbox\big(B_1^G(e)\big)-\dimbox\big(B_1^{U^+}(e)\big)-
\dimbox\big(B_1^{U^-}(e)\big)$.
\end{proposition}

The proof of Proposition~\ref{pro: box dimension}
follows closely
\cite[Prop.~9.3]{Einsiedler-Katok-Lindenstrauss} which gives a
similar statement for $\SL(k,\RR)/\SL(k,\ZZ)$ and $k\geq 3$.

We recall the notion of topological entropy for a
continuous transformation $\beta:Y\rightarrow Y$ on a compact
metric space $(Y,d)$. For $\epsilon>0$ and a positive integer $n$
let $s_{\epsilon,n}(\beta)$  be the maximal number of points
$y_1,y_2,\ldots$ such that for any $i\neq j$ there exists some
$0\leq k<n$ with $d(\beta^ky_i,\beta^ky_j)>\epsilon$. Then the
{\em topological entropy} is defined by
\[
   \htop(\beta)=\lim_{\epsilon\rightarrow0}
    \limsup_{n\rightarrow\infty}\frac{s_{\epsilon,n}(\beta)}{n}.
\]

\begin{lemma}\label{lemma topological entropy}
   Suppose Proposition \ref{pro: box dimension} were not true. Then
   the restriction $\beta$ of $\alpha^{(\bt,\bn)}$ to the
   $\alpha|_C$-invariant compact set $Y_C$
   has positive topological entropy.
\end{lemma}

\begin{proof}
  It is easy to check that $Y_C$ is compact and $\alpha^{(\bt,\bn)}$-invariant
  for any $(\bt,\bn)\in C$. Also note, that there exists $\delta>0$
  such that for any $y\in Y$ the map $g\in B_{2\delta}^G(e)\mapsto gy\in X$
  is an isometry. For simplicity write $a=\alpha^{(\bt,\bn)}\in G$.
  The Lie algebra $\mathfrak u^+$ of $U^+$ is
  generated by the eigenvectors of $\Ad_a$ to eigenvalues of absolute 
value bigger than
  one. Similarly, $\mathfrak u^-$ is the Lie algebra to $U^-$, and let
  $\mathfrak u^0$ be the linear hull of the eigenspaces to all 
eigenvalues of $\Ad_a$
  with absolute value one.
  Therefore, there exists $\lambda>1$ and $c>0$ such that
  $\|\Ad_a^j w\|\geq c \lambda^j\|w\|$ for all $w\in\mathfrak u^+$, 
where $\|\cdot\|$
  is some fixed norm on the Lie algebra $\mathfrak u^+$ -- for more details
  see \cite[Sect.~4]{Einsiedler-Katok-nonsplit}.
  Also note that $\exp (\Ad_a u)=a\exp(u)a^{-1}$.

  Recall the assumption
  that the right invariant metric $d(\cdot,\cdot)$ on $G$ makes $\exp$ 
a bi-Lipshitz map.
  We conclude that there exists some $\epsilon>0$
  such that for $u,v\in U^+$ with $\lambda^{-n}/2<d(u,v)<\delta$
  there exists $j$ with $0\leq j< n$  and
  $\epsilon<d(a^ju,a^jv)=d(a^jua^{-j},a^jva^{-j})<\delta$.
  The same holds similarly for $u,v\in U^-$.

  We first consider the statement about $E_{B,x}=E_{B,x}(Y)$
  in Proposition~\ref{pro: box dimension} where
  $B=B_{\delta/2}^{U^+}(e)\subset U^+$ and $x\in X$.
  Then the above shows that for any two points $u,v\in E_{B,x}$
  with $d(u,v)>\lambda^{-n}$ there exists a $j$ with $0\leq j<n$
  such that $d(\beta^j ux,\beta^j vx)>\epsilon$.
  Note that by definition of $E_B$ we have $\beta^j ux,\beta^j vx\in Y$
  for all $j\geq 0$.
  Therefore $s_{\epsilon,n}(\beta)\geq s_{\lambda^{-n}}(E_{B,x})$.
  With this it  is easy to derive the lemma in this case. In fact,
   \begin{multline*}
    \htop(\beta)\geq\limsup_{n\rightarrow\infty}\frac{\log
  s_{\epsilon,n}(\beta)}{n}
   \geq(\log\lambda)\limsup_{n\rightarrow\infty}
           \frac{\log s_{\lambda^{-n}}\big(t_V(E_{B,x},x)\big)}{|\log\lambda^{-n}|}=\\
  (\log\lambda) \limsup_{\delta\rightarrow 0}\frac{\log 
s_\delta\big(t_V(E_{B,x},x)\big)}{|\log\delta|}
      =(\log\lambda)\dimbox\big(t_V(E_{B,x},x)\big)>0.
   \end{multline*}
  In general $B\subset U^+$ is assumed to be compact, and therefore
  can be covered by balls as above. This shows the lemma in the case
  of the first statement of Proposition \ref{pro: box dimension}.

  So assume now $-(\bt,\bn)\in C$. Let $\eta>0$ and let
  $B=B_\eta^{U^+}(g^+)B_\eta^{U^-}(g^-)Tg^0$, where
  $T=\exp \bigl(B^{\mathfrak u^0}_\eta(0)\bigr)$
  is the exponential image of a ball in $\mathfrak u^0$.
  Note that $B\subset G$ is open, and that $B\subset B_\delta(g^+g^-g^0)$
  for small enough $\eta$. Moreover, the metric restricted to $B$
  is Lipshitz equivalent to the product metric on
  $B_\eta^{U^+}(g^+)\times B_\eta^{U^-}(g^-)\times (Tg^0)$.
  If the last statement of Proposition \ref{pro: box dimension} fails,
  then the image $I_{B,x}$ of the projection of $E_{B,x}=\{u\in B: 
ux\in Y_C\}$ to
  $B_\eta^{U^+}(g^+)\times B_\eta^{U^-}(g^-)$ has positive box dimension.
  We show as before that this implies that $\htop(\beta)>0$.
If $u,v\in B$ and $u=u^+u^-u^0$ (resp., $v=v^+v^-v^0$) are
  such that $d(u^+u^-,v^+v^-)>\lambda^{-n}$, then either
  $d(u^+,v^+)>\lambda^{-n}/2$ or $d(u^-,v^-)>\lambda^{-n}/2$.
  In the former case we find $j$ with $0\leq j<n$,
  $\epsilon<d(a^ju^+a^{-j},a^jv^+a^{-j})<\delta$, and 
$d(a^ju^-u^0a^{-j},a^jv^-v^0a^{-j})<c\eta$
  for some absolute constant $c$. For $\epsilon'= \epsilon/2$ this implies
  for small enough $\eta$ that
  $\epsilon'<d(a^j u a^{-j},a^jva^{-j})<2\delta$. The second case
  leads to a similar estimate for some $j$ with $-n<j\leq 0$.
  Using the invariance of $Y_C$ under $\alpha^{-(\bt,\bn)}$
  this shows that $s_{\epsilon',2n-1}(\beta)\geq s_{\lambda^{-n}}(I_{B,x})$.
  As above this gives the desired result.
\end{proof}

We can link Lemma \ref{lemma topological entropy} with the
measure-theoretic entropy by the
varational principle (see \cite[Thm.~8.6]{Walters-82}): For any
continuous map $\beta:Y\rightarrow Y$ of a compact metric space $Y
$, its topological entropy is the supremum
\begin{equation}\label{varational principal}
    \htop(\beta)=\sup_\nu\h_\nu(\beta)
\end{equation}
over all measure-theoretic entropies with respect to
$\beta$-invariant probability measures $\nu $ on $Y $.

\begin{proof}[Proof of Proposition \ref{pro: box dimension}]
By Lemma \ref{lemma topological entropy}
it suffices to assume that $\htop(\beta)>0$
and to conclude a contradiction.

By the varational principle \eqref{varational principal} there
exists a $\beta$-invariant measure $\nu$ (supported on the compact
set $Y_C $ where $\beta$ is defined) such that $\h_\nu(\beta)$ is
as close to $\htop(\beta)$ as we want.  All we need is
$\h_\nu(\beta)>0$. Since the entropy of $\beta $ with respect to a
non-ergodic measure $\nu$ can be calculated as an integral of the
entropies of $\beta $ with respect to the ergodic components of
$\nu $, we can assume that $\nu$ is $\beta $-ergodic.

There is no reason for $\nu $ to be $\alpha$-invariant, so we need
to go through the following averaging procedure. Let $\lambda $
denote the Haar measure on $\RR^k\times\ZZ^\ell$ and let
$Q_N=[-N,N]^k\times\{-N,\ldots,N\}^\ell$ be a large box in
$\RR^k\times\ZZ^l$. Define the sequence of measures
\begin{equation}\label{definition of nu N}
   \nu_N=\frac{1}{\lambda(Q_N\cap C)}\int_{Q_N\cap C}
    (\alpha^{(\bs,\bm)})_*\nu\operatorname{d}\!\lambda
\end{equation}
on $Y_C $.  Then the definition of the cone implies that for large $N $ the
measure $(\alpha^{(\bs,\bm)})_*\nu_N$ is close (in strong or weak$^*$
topology) to $\nu_N $. Therefore, every weak$^*$ accumulation
point $\bar\nu$ of the sequence $\nu_N$ will be an $\alpha
$-invariant probability measure on $b_C $.

In many different settings it is well known (see
   \cite[Thm.~8.2]{Walters-82} and \cite[Thm.~4.1]{Newhouse-89}) that the
measure-theoretic entropy depends upper semi-\-con\-tin\-u\-ously
on the measure (using the weak$^*$ topology). These theorems do
not quite apply to the proof here, but it is easy to adapt the
proof of \cite[Cor.~10.3]{Einsiedler-Katok-Lindenstrauss} to
$\beta$. By definition of $\nu_N$ in \eqref{definition of nu N}
the ergodic components of $\nu_N$ with respect to $\beta$ are all
of the form $(\alpha^{(t,n)})_*\nu$. Therefore,
$\h_{\nu_N}(\beta)=\h_\nu(\beta)$ for all $N$, and upper
semi-continuity implies that
\[
   \h_{\bar\nu}(\beta)\geq\h_\nu(\beta)>0.
\]
Again because entropy can be expressed as an integral -- this time
over all ergodic components $\mu$ of $\bar\nu$ with respect to
$\alpha$ -- there exists an $\alpha$-invariant and ergodic measure
$\mu$ that is supported by $Y$ with $\h_\mu(\alpha^{(\bt,\bn)})>0$.
This is a contradiction to the assumption of Proposition \ref{pro: 
box dimension}.
\end{proof}

\section{Conclusion}

\label{concl}

\subsection{Completion of the proof of Theorem \ref{theorem: Main}}
\label{proof}
\begin{proof}[Proof of Theorem \ref{theorem: Main}]
  We first consider the case $v\in\ZZ_p$ and define
  \[
   E_v=\bigl\{u\in\RR:\liminf_{q\rightarrow\infty, 
q_0\rightarrow\infty} |q|\cdot|qu-q_0|\cdot|qv-q_0|>0\bigr\}.
  \]
  Here Proposition \ref{prop: dynamical g} shows $u\in E$ if and only 
if $\alpha^C x_{u,v}$ is unbounded, i.e.~ $\alpha^Cx_{u,v}\subset 
K_\delta$ for some $\delta>0$,
where $K_\delta$ are the compact subsets defined in \S \ref{sec: 
Mahler}. Therefore,
  \[
   E_v=\bigcup_{\delta>0}E_{v}(\delta)
  \]
  where
  \[
   E_v(\delta)=\left\{u\in[-1/\delta,1/\delta]: 
\alpha^C\left(\begin{pmatrix} 1 &\\ u&1\end{pmatrix}, \begin{pmatrix} 
1 &\\ &1\end{pmatrix}\right)x_{0,v}\subset K_\delta\right\}.
  \]
  By Corollary \ref{corollary}  $K_\delta\neq X$ does not support any 
$\alpha$-invariant and ergodic probability measure $\mu$ with 
$\h_\mu(\alpha^{(1,0)})>0$. It is easy to check that
  \[
   U=\left\{\left(\begin{pmatrix} 1 &\\ u&1\end{pmatrix}, 
\begin{pmatrix} 1 &\\ &1\end{pmatrix}\right):u\in\RR\right\}
  \]
  is the horospherical unstable subgroup corresponding to 
$\alpha^{(1,0)}$. Therefore, $E_v(\delta)$ has box dimension zero by 
Proposition \ref{pro: box dimension}, and so $E_v$ is a countable 
union of such sets as claimed in Theorem \ref{theorem: Main}. The 
fact that $E_v$ has Hausdorff dimension zero follows from this: 
Hausdorff dimension is always less or equal to box dimension, and the 
Hausdorff dimension of a countable union is the supremum of the 
Hausdorff dimensions of the original sets.

Now let $v\in \QQ_p$ and fix some $n\geq 0$ with $p^n v\in \ZZ_p$. 
Using the definition only, it is straightforward to check that $p^n 
E_v\subset E_{p^n v}$. Since the latter satisfies the conclusion of 
the theorem, so does $E_v$.
\end{proof}

\subsection{Furstenberg's theorem and Theorem 
\ref{t:Furstenberg}}\label{s:Furstenberg}

Let $\TT=\RR/\ZZ$ and denote by $\times p$ the map $x\mapsto px$ (modulo one).
Recall that a closed subset $F\subset\TT$ is $\times p$-invariant if
$\times p(F)\subset F$.

\begin{theorem}\cite[Thm.~IV.1]{Furstenberg67}\label{thmF}
  Let $p_1,p_2$ be two distinct primes. A closed subset $F\subset\TT$ 
that is $\times p_1,\times p_2$-invariant
is either a finite set of rational points or equal to $\TT$.
\end{theorem}

\begin{proof}[Proof of Theorem \ref{t:Furstenberg}]
For two primes $p_1, p_2$, and $\delta>0$, consider  the set
\[
  F_\delta=\bigl\{u\in\RR: 
\inf_{q\in\NN}|q|\cdot|q|_{p_1}\cdot|q|_{p_2}\cdot\langle 
qu\rangle\geq\delta\bigr\}/\ZZ.
\]
 From the definition it is clear that $F_\delta$ is closed. We claim it
is also $\times p_1,\times p_2$-invariant. Indeed,  let $u\in 
F_\delta$ and $q\in\NN$,
then
\[
|q|\cdot|q|_{p_1}\cdot|q|_{p_2}\cdot\langle qp_1u\rangle= 
|qp_1|\cdot|qp_1|_{p_1}\cdot|qp_1|_{p_2}\cdot\langle 
qp_1u\rangle\geq\delta.
\]
Since $q\in\NN$ was arbitrary, $p_1 u\in F_\delta$.
Clearly, $F_\delta$ does not contain any rational numbers.
By Theorem \ref{thmF} $F_\delta$ must be empty, which shows that for 
every $u\in\RR$
there exists $q\in\NN$ with
\[
  |q|\cdot|q|_{p_1}\cdot|q|_{p_2}\cdot\langle qu\rangle<\delta.
\]
If $u$ is rational there is nothing to prove. Otherwise,
it is clear that for small $\delta>0$ the natural number $q$ as above
needs to be rather big. Therefore, \equ{e:MTmF} is satisfied
and the theorem follows.
\end{proof}

\subsection{$\cD$-adic valuations}\label{s:D} Following \cite{Mathan-Teulie},
let us consider an increasing sequence $\cD = ( r_n)_{n\in\ZZ_+}$
of positive integers with $r_0 = 1$ and $r_{n+1}
\in r_n\ZZ$ for each $n$, and define
  the $\cD$-adic pseudo-valuation
$$
|q|_{\cD} \df \inf\{1/r_n : q\in r_n\ZZ\}\,.
$$
Clearly $
|\cdot|_{\cD} = |\cdot|_p$ for $\cD = \{p^n\}$, in which case $
|\cdot|_{\cD}$ is multiplicative.
De~Mathan and Teuli\'e showed that 
whether equality
\begin{equation}\label{e:MTma}
  \liminf_{q\rightarrow\infty }|q|\cdot|q|_{\cD}\cdot\langle qu\rangle=0
\end{equation}
is equivalent to the digits in the
continued fraction expansions for $r_nu$ not having a uniform bound. 
They also proved
that \equ{e:MTma}
holds when $u$ is quadratic irrational and
\begin{equation}\label{e:bdd}
\text{ratios }r_{n+1}
/r_n\text{ are uniformly bounded,}
\end{equation}
and asked whether or not it holds for all real numbers $u$.

It seems tempting to attack this problem using the methods of the 
present paper.
Unfortunately for general $\cD$, even satisfying \equ{e:bdd}, we do 
not know how
to translate \equ{e:MTma} into the language of group actions. However 
the special case
$\cD = \{a^n\}$, where $a$ is not necessarily prime, happens to be amenable
to our technique and we can assert that for this choice of $\cD$
the set of exceptions to \equ{e:MTma} has \hd\ zero.
Further, if we denote by $S$ the set $\{p_1,\dots,p_\ell\}$ of primes 
dividing $a$,
the completion of $\QQ$
with respect to $
|\cdot|_{\cD}$ is bi-Lipschitz equivalent to $\QQ_S \df \prod_{i}\QQ_{p_i}$
(where
the metric on $\QQ_{p_i}$ is scaled by the exponent with which $p_i$ 
occurs in $a$).
Then it can be shown that
for any  $v\in\QQ_S$,
the set of $u\in\RR$
which do not satisfy
\begin{equation}\label{g-conjecture-a}
  \liminf_{
q\rightarrow\infty,\,q_0\rightarrow\infty
}  |q|\cdot|qu-q_0|\cdot|qv-q_0|_\cD=0
\end{equation}
  is a
  countable union of sets with box dimension zero.

The proof  uses dynamics on $X=G/\Gamma$, where
$$
  G = \SL(2,\RR)\times\SL(2,\QQ_S)= \SL(2,\RR)\times\prod_i\SL(2,\QQ_{p_i})\,,
$$
  and   $\Gamma$  is
$\SL\big(2,\ZZ[\tfrac{1}{p_1},\dots,\tfrac{1}{p_\ell}]\big)$
  diagonally embedded in $G$.   The action considered is left multiplication  by
$
  \left(\begin{pmatrix}e^{-t}&\\&e^t\end{pmatrix},
\begin{pmatrix}a^{n}&\\&a^{-n}\end{pmatrix}\right)\mbox{ for 
}(t,n)\in\RR\times\ZZ
$, where $\begin{pmatrix}a^{n}&\\&a^{-n}\end{pmatrix}$ is understood 
to be diagonally
imbedded in $\SL(2,\QQ_S)$.
One can prove  that  $u \ne 0$ (resp., $u \ne 0$ and 
$v\in\prod_i\SL(2,\ZZ_{p_i})$)
  fail
\eqref{e:MTma} (resp.,  \eqref{g-conjecture-a}) if and only if certain elements
of $X$ as above have bounded orbits relative to certain cones in 
$\RR\times\ZZ$.
Then one applies Theorem \ref{theorem: Elon} with $L = \SL(2,\QQ_S)$.

\subsection{Extending Theorem \ref{theorem: Main}}\label{forthcoming}

Theorem \ref{theorem: Main} can be strengthened in another direction.
Namely, the following is true:

\begin{theorem}\label{theorem: More} The set of $\,(u,v)\in\RR\times\QQ_p$
which do not satisfy
\begin{equation}\label{g-conjecture-q}
  \liminf_{
q\rightarrow\infty,\,q_0\rightarrow\infty
}  \,|q|\cdot|qu-q_0|\cdot|qv-q_0|_p=0
\end{equation}
  is at most a
  countable union of sets with box dimension zero.
\end{theorem}

For the above a stronger version of Theorem \ref{theorem: Elon} is 
needed, namely that the Haar measure $m_X$ is the only 
$\alpha$-invariant and ergodic probability measure with 
$\h_\mu(\alpha^{(t,n)})>0$ for some $(t,n)\in\RR\times\ZZ$. The  case 
$(t,n)=(0,1)$ opposite to Theorem \ref{theorem: Elon} where 
$(t,n)=(1,0)$ can be proven with the same method
as the main theorem of  \cite{Lindenstrauss-Quantum}. However, we claim that
\[
  \h_\mu(\alpha^{(t,n)})=\h_\mu(\alpha^{(t,0)})+\h_\mu(\alpha^{(0,n)}),
\]
and so the above two cases give the complete answer for measures with 
positive entropy.
This additivity of entropy follows from 
\cite{Einsiedler-Katok-nonsplit}, since
for a general $(t,n)$, the stable horospherical subgroup $U$ is the product of
its real subgroup $U_\infty$ and its $p$-adic subgroup $U_p$, and  the Lyapunov
weights that describe the contraction (resp., expansion) rate for $U_\infty$
(resp., $U_p$) are linearly independent (indeed, $\alpha^{(1,0)}$ commutes with
the $p$-adic group $U_p$ but not with $U_\infty$, and similarly for 
$\alpha^{(0,1)}$).
By \cite[Thm.~8.4]{Einsiedler-Katok-nonsplit} this shows that the 
conditional measures
for $U$ are the product measures of the conditional measures for 
$U_\infty$ and $U_p$.
By \cite[Prop.~9.4]{Einsiedler-Katok-nonsplit} these conditional 
measures determine
entropy which implies the above.

A more general measure rigidity theorem, covering this as well as 
other more general
situations, is an ongoing joint work of E.~Lindenstrauss and the 
first named author.

\def\cprime{$'$}
\providecommand{\bysame}{\leavevmode\hbox to3em{\hrulefill}\thinspace}
\providecommand{\MR}{\relax\ifhmode\unskip\space\fi MR }
\providecommand{\MRhref}[2]{%
   \href{http://www.ams.org/mathscinet-getitem?mr=#1}{#2}
}
\providecommand{\href}[2]{#2}

\end{document}